\documentclass[11pt]{article}

\usepackage[latin1]{inputenc}

\usepackage{amsfonts,amsmath,amssymb,amsthm,graphicx,epsfig,float}
\usepackage[T1]{fontenc}

\usepackage[english,francais]{babel}

{
  \newtheorem{theoreme}{Th\'eor\`eme}
  \newtheorem*{theoreme*}{Th\'eor\`eme}
  \newtheorem{lemme}[theoreme]{Lemme}

  \newtheorem{proposition}[theoreme]{Proposition}
\newtheorem*{corollaire*}{Corollaire}
\newtheorem*{proposition*}{Proposition}
\theoremstyle{remark}
  \newtheorem*{remarque*}{Remarque}
}

\newcounter{ex}

\newenvironment{rem*}{
  \noindent\textbf{Remarque. }}{}

\newcommand{\Cc}{\mathbb{C}}

\newcommand{\Pp}{\mathbb{P}}

\title{{\bf Sur les automorphismes réguliers de $\Cc^k$}}
\author{Henry de Thélin}
\date{}

\begin{document}
\maketitle


\def\figurename{{Fig.}}%
\def\proofname{Preuve}
\def\contentsname{Sommaire}%

\begin{abstract}

Nous montrons l'unicité de la mesure d'entropie maximale pour les automorphismes réguliers de $\Cc^k$.

\end{abstract}

\selectlanguage{english}
\begin{center}
{\bf{ }}
\end{center}

\begin{abstract}

We show the uniqueness for the measure of maximal entropy for regular automorphisms of $\Cc^k$.

\end{abstract}

\selectlanguage{francais}

Mots-clefs: dynamique complexe, courants, entropie.

Classification: 37Fxx, 32H50, 37A35, 37Dxx.

\section*{{\bf Introduction}}
\par

Soit $f: X \rightarrow X$ une application holomorphe ou méromorphe sur une variété complexe compacte $X$. Quand on sait construire une mesure invariante d'entropie maximale pour $f$, une question naturelle est de savoir si $f$ en possède d'autres.

Lorsque $f$ est une application de Hénon de $\Cc^2$ (voir \cite{BLS}), un Hénon-Like (voir \cite{Du}), ou un endomorphisme holomorphe de $\Pp^k(\Cc)$ (voir \cite{Ly} et \cite{Ma} pour la dimension $1$ et \cite{BD} pour la dimension supérieure), $f$ possède une unique mesure d'entropie maximale.

Dans cet article, nous nous intéressons à cette question pour les automorphismes réguliers de $\Cc^k$ (voir \cite{S} ou le paragraphe \ref{rappels} pour la définition).

Pour $f$ un automorphisme régulier de $\Cc^k$, N. Sibony a construit dans \cite{S} la mesure d'équilibre $\mu$ et celle-ci est d'entropie maximale (voir \cite{G}). L'objectif de cet article est alors de démontrer le théorème suivant:

\begin{theoreme}

La mesure $\mu$ est l'unique mesure d'entropie maximale de $\Cc^k$ pour $f$.

\end{theoreme}

Pour démontrer ce théorème, nous suivrons l'approche utilisée par E. Bedford, M. Lyubich et J. Smillie dans le cas des applications de Hénon (voir \cite{BLS}). Il s'agira de changer un certain nombre d'arguments qui sont propres à la dimension $2$.

\section{{\bf Rappels}}{\label{rappels}}

Dans ce paragraphe, nous faisons des rappels sur les
automorphismes réguliers de $\Cc^k$ et nous donnons quelques
propriétés qui nous serviront notre démonstration du théorème.

\subsection{{\bf Automorphismes réguliers de $\Cc^k$ }}{\label{automorphismes}}

Soit $f$ un automorphisme polynomial de $\Cc^k$. On peut étendre $f$
en une application birationnelle de $\Pp^k(\Cc)$ et nous noterons
encore $f$ cette extension. Soient $I_{+}$ et $I_{-}$ les ensembles
  d'indétermination respectivement de $f$ et $f^{-1}$. L'automorphisme
  $f$ est dit régulier lorsque $I_{+} \cap I_{-} = \emptyset$ (voir \cite{S}).

Dans toute la suite, nous considèrerons $f$ un automorphisme régulier
de $\Cc^k$. Voici quelques unes de ses propriétés (voir \cite{S} et \cite{DS5}).

Il existe un entier $s$ tel que $\mbox{dim } I_{+}=k-s-1$ et
$\mbox{dim } I_{-}=s-1$. Ces ensembles sont contenus dans l'hyperplan à
l'infini $L_{\infty}$ et on a $f(L_{\infty} \setminus
I_{+})=f(I_{-})=I_{-}$ et $f^{-1}(L_{\infty} \setminus
I_{-})=f^{-1}(I_{+})=I_{+}$. Si $d_{\pm}$ désignent les degrés
algébriques de $f$
et $f^{-1}$, alors $d_{+}^s=d_{-}^{k-s}$.

Rappelons maintenant la définition des degrés dynamiques $d_q$ de $f$
(voir \cite{RS}). On pose $\delta_q(f):= \int_{\Pp^k(\Cc)} f^{*} (\omega^q)
\wedge \omega^{k-q}$ pour $q=0 \mbox{,} \dots \mbox{,} k$ et on a:

$$d_q:= \lim_{n \rightarrow \infty} (\delta_q(f^n))^{1/n}.$$

Les automorphismes réguliers de $\Cc^k$ sont des cas particuliers des
applications birationnelles régulières de T.-C. Dinh et N. Sibony
(voir \cite{DS4}). Il résulte de leur article
\cite{DS4} que $d_q= d_{+}^q$ pour $q=0 \mbox{,} \dots \mbox{,} s$ et $d_q=
d_{-}^{k-q}$ si $q=s \mbox{,} \dots \mbox{,} k$. 

Soient $K_{+}$ (respectivement $K_{-}$) l'ensemble des points $z$ de
$\Cc^k$ pour lesquels l'orbite $(f^n(z))_{n \geq 0}$ (respectivement
$(f^{-n}(z))_{n \geq 0}$) est bornée. On a $\overline{K_{\pm}}=
K_{\pm} \cup I_{\pm}$. L'ensemble $I_{+}$ est attirant pour $f^{-1}$: cela
signifie qu'il existe un voisinage $V_{+}$ de $I_{+}$ tel que
$f^{-1}(V_{+}) \Subset V_{+}$ et $\cap_{n \geq 0}
f^{-n}(V_{+})=I_{+}$. De même, l'ensemble $I_{-}$ est attirant pour
$f$ et on notera $V_{-}$ un voisinage de $I_{-}$ disjoint de $V_{+}$,
avec $ f(V_{-}) \Subset V_{-}$ et $\cap_{n \geq 0}
f^{n}(V_{-})=I_{-}$. Il découle de ces propriétés que le bassin
d'attraction pour $f^{-1}$ de $I_{+}$ est égal à $\Pp^k(\Cc) \setminus
\overline{K_{-}}$ et que le bassin d'attraction pour $f$ de $I_{-}$
est égal à  $\Pp^k(\Cc) \setminus \overline{K_{+}}$.

La suite de $(1,1)$ courant $d_{\pm}^n (f^{\pm n})^{*} \omega$
converge vers le courant de Green $T_{\pm}$. Rappelons brièvement une
démonstration (celle de \cite{GS}) de cette convergence car elle nous
sera utile dans la suite.

On a $\frac{f^{*} \omega}{d_{+}}= \omega + dd^c u$ où $u$ est une
fonction quasi-psh qui est lisse en dehors de $I_{+}$ et qui peut être supposée négative. En itérant cette relation, on a:

$$T_{n,+}=\frac{(f^n)^{*} \omega}{d_{+}^n}= \omega + dd^c \sum_{i=0}^{n-1}
\frac{u \circ f^i}{d_{+}^i}.$$

La suite $v_n=\sum_{i=0}^{n-1} \frac{u \circ f^i}{d_{+}^i}$ est bornée en dehors de $V_{+}$. Elle décroît donc vers une fonction
quasi-psh $v_{\infty}$ qui est continue en dehors
de $I_{+}$. De plus la suite $\omega + dd^c v_n$ converge au sens des
courants vers $T_{+}:=\omega + dd^c v_{\infty}$.

Les courants $T_{\pm}$ vérifient $f^{*} T_{+} = d_{+} T_{+}$ et $f_{*}
T_{-} = d_{-} T_{-}$. Les puissances $T_{+}^s$ et $T_{-}^{k-s}$ sont
bien définies. T.-C. Dinh et N. Sibony ont démontré dans \cite{DS5}
que $T_{+}^s$ est l'unique courant positif fermé de bidegré $(s,s)$ et
de masse $1$ qui a son support dans $\overline{K_{+}}$. De même,
$T_{-}^{k-s}$ est l'unique courant positif fermé de bidegré
$(k-s,k-s)$ et de masse $1$ qui a son support dans
$\overline{K_{-}}$. 

La mesure $\mu= T_{+}^s \wedge T_{-}^{k-s}$ est bien définie et son
support est dans le bord de $K:= K_{+} \cap K_{-}$. Cette mesure est une probabilité invariante, mélangeante et d'entropie maximale $\log d_s= \log d_{+}^s$ (voir
\cite{G}). 

Dans le paragraphe suivant, nous allons donner quelques propriétés techniques dont on se servira dans notre démonstration.

\subsection{{\bf Quelques propriétés}}{\label{propriete}}

Dans ce paragraphe on considère un ouvert $\Omega$ inclus dans $(\Pp^k(\Cc) \setminus V_{+}) \cap
(\Pp^k(\Cc) \setminus V_{-})$ et $\xi$ un sous-ensemble analytique
lisse de $\Omega$ de dimension pure $s$ qui se prolonge un petit peu. Comme $f^{-1}(V_{+}) \Subset V_{+}$ et $ f(V_{-})
\Subset V_{-}$, $f^n(\Omega)$ vit dans $\Pp^k(\Cc) \setminus V_{+}$
et $f^{-n}(\Omega)$ est dans $\Pp^k(\Cc) \setminus V_{-}$. En
particulier $f^n$ et $f^{-n}$ sont holomorphes sur $\Omega$ pour $n \geq 0$.

$T_{+}$ possède un potentiel qui est continu sur $\Omega$, on peut
donc définir par récurrence les wedges $[\xi] \wedge T_{+}^p$ pour $p$
compris entre $0$ et $s$ (voir
\cite{BT}). Nous allons donner quelques propriétés techniques sur ces wedges dont on se servira dans nos démonstrations.

\begin{lemme}

Soit $0 \leq \psi \leq 1$ une fonction $C^{\infty}$ à support compact
dans $\Omega$. Alors, on a: 

$$I=\int \psi [\xi] \wedge \omega^l \wedge \frac{(f^{i_1})^{*} \omega}{d_{+}^{i_1}} \wedge \dots \wedge \frac{(f^{i_p})^{*} \omega}{d_{+}^{i_p}} \wedge \frac{(f^n)^{*}
T_{m,+}^{s-l-p} }{d_{+}^{(s-l-p)n}} \leq K(\xi, \psi)$$

pour tous entiers $l$ et $p$ compris entre $0$ et $s$ avec $0 \leq l+p
\leq s$ et tous entiers $i_1, \dots , i_p$, $m$ et $n$. De plus on a la même majoration si on remplace $T_{m,+}$ par $T_{+}$.

\end{lemme}

\begin{proof}

On fait la preuve pour $T_{m,+}$ car c'est la même pour $T_{+}$.

On va commencer par remplacer les $\frac{(f^n)^{*} T_{m,+}}{d_{+}^n}$ par des $\omega$ ou des
$\frac{(f^n)^{*} \omega}{d_{+}^n}$.

On a $T_{m,+}= \omega + dd^c v_m$ (voir les rappels précédents pour
les notations). En particulier,

\begin{equation*}
\begin{split}
&I = \int \psi [\xi] \wedge \omega^{l} \wedge \frac{(f^{i_1})^{*} \omega}{d_{+}^{i_1}} \wedge \dots \wedge \frac{(f^{i_p})^{*} \omega}{d_{+}^{i_p}} \wedge \frac{(f^n)^{*} \omega}{d_{+}^n}  \wedge \frac{(f^n)^{*} T_{m,+}^{s-l-p-1}}{d_{+}^{(s-l-p-1)n}}\\
&+ \int \frac{v_m \circ f^n}{d_{+}^n} dd^c \psi \wedge [\xi] \wedge \omega^l \wedge \frac{(f^{i_1})^{*} \omega}{d_{+}^{i_1}} \wedge \dots \wedge \frac{(f^{i_p})^{*} \omega}{d_{+}^{i_p}}  \wedge \frac{(f^n)^{*} T_{m,+}^{s-l-p-1}}{d_{+}^{(s-l-p-1)n}}.
\end{split}
\end{equation*}

Mais comme $f^r(\Omega)$ est dans $\Pp^k(\Cc) \setminus V_{+}$ pour
tout $r \geq 0$ et
que $u$ est bornée sur cet ensemble, on en déduit que $|v_m \circ f^n|
\leq K $ dans $\Omega$. Cela implique que 

\begin{equation*}
\begin{split}
&\left| \int \frac{v_m \circ f^n}{d_{+}^n} dd^c \psi \wedge [\xi] \wedge \omega^l \wedge \frac{(f^{i_1})^{*} \omega}{d_{+}^{i_1}} \wedge \dots \wedge \frac{(f^{i_p})^{*} \omega}{d_{+}^{i_p}}  \wedge \frac{(f^n)^{*} T_{m,+}^{s-l-p-1}}{d_{+}^{(s-l-p-1)n}} \right| \leq\\
&\frac{K(\psi)}{d_{+}^n}  \int \psi_1  [\xi] \wedge \omega^{l+1} \wedge \frac{(f^{i_1})^{*} \omega}{d_{+}^{i_1}} \wedge \dots \wedge \frac{(f^{i_p})^{*} \omega}{d_{+}^{i_p}}  \wedge \frac{(f^n)^{*} T_{m,+}^{s-l-p-1}}{d_{+}^{(s-l-p-1)n}},
\end{split}
\end{equation*}

où $\psi_1$ est une fonction $C^{\infty}$ à support compact dans
$\Omega$ et comprise entre $0$ et $1$. On obtient alors

\begin{equation*}
\begin{split}
&I \leq \int \psi [\xi] \wedge \omega^{l} \wedge \frac{(f^{i_1})^{*} \omega}{d_{+}^{i_1}} \wedge \dots \wedge \frac{(f^{i_p})^{*} \omega}{d_{+}^{i_p}} \wedge \frac{(f^n)^{*} \omega}{d_{+}^n}  \wedge \frac{(f^n)^{*} T_{m,+}^{s-l-p-1}}{d_{+}^{(s-l-p-1)n}} +\\
&K(\psi)  \int \psi_1  [\xi] \wedge \omega^{l+1} \wedge \frac{(f^{i_1})^{*} \omega}{d_{+}^{i_1}} \wedge \dots \wedge \frac{(f^{i_p})^{*} \omega}{d_{+}^{i_p}}  \wedge \frac{(f^n)^{*} T_{m,+}^{s-l-p-1}}{d_{+}^{(s-l-p-1)n}},
\end{split}
\end{equation*}

et on a ainsi éliminé un $T_{m,+}$. En recommençant ce procédé avec
les deux intégrales ci-dessus, on prouve (quitte à renommer la
constante $K(\psi)$):

$$I \leq K(\psi) \sum_{l'=l}^{s-p} \int \phi_{l'} [\xi] \wedge
  \omega^{l'} \wedge  \frac{(f^{i_1})^{*} \omega}{d_{+}^{i_1}} \wedge \dots \wedge \frac{(f^{i_p})^{*} \omega}{d_{+}^{i_p}} \wedge \frac{(f^n)^{*} \omega^{s-p-l'}}{d_{+}^{(s-p-l')n}},$$

où les $\phi_{l'} $ sont $C^{\infty}$ comprise entre $0$ et $1$ et à
support compact dans $\Omega$.

Maintenant, en utilisant la même méthode, on va enlever les $\frac{(f^i)^{*} \omega}{d_{+}^i}$.

Notons $J= \int \phi_{l'} [\xi] \wedge  \omega^{l'} \wedge  \frac{(f^{i_1})^{*} \omega}{d_{+}^{i_1}} \wedge \dots \wedge \frac{(f^{i_p})^{*} \omega}{d_{+}^{i_p}} \wedge \frac{(f^n)^{*} \omega^{s-p-l'}}{d_{+}^{(s-p-l')n}}$. On a, si $l' < s$, par exemple

\begin{equation*}
\begin{split}
&J= \int \phi_{l'} [\xi] \wedge \omega^{l'+1} \wedge   \frac{(f^{i_2})^{*} \omega}{d_{+}^{i_2}} \wedge \dots \wedge \frac{(f^{i_p})^{*} \omega}{d_{+}^{i_p}} \wedge \frac{(f^n)^{*} \omega^{s-p-l'}}{d_{+}^{(s-p-l')n}} +\\
&\int v_{i_1} dd^c  \phi_{l'} \wedge [\xi] \wedge
  \omega^{l'} \wedge  \frac{(f^{i_2})^{*} \omega}{d_{+}^{i_2}} \wedge \dots \wedge \frac{(f^{i_p})^{*} \omega}{d_{+}^{i_p}} \wedge \frac{(f^n)^{*} \omega^{s-p-l'}}{d_{+}^{(s-p-l')n}} .
\end{split}
\end{equation*}

Comme précédemment, 

\begin{equation*}
\begin{split}
& \left| \int v_{i_1} dd^c  \phi_{l'} \wedge [\xi] \wedge
  \omega^{l'} \wedge  \frac{(f^{i_2})^{*} \omega}{d_{+}^{i_2}} \wedge \dots \wedge \frac{(f^{i_p})^{*} \omega}{d_{+}^{i_p}} \wedge \frac{(f^n)^{*} \omega^{s-p-l'}}{d_{+}^{(s-p-l')n}}  \right| \leq\\
 & K(\psi) \int \phi_{l',1}  [\xi] \wedge
  \omega^{l'+1} \wedge    \frac{(f^{i_2})^{*} \omega}{d_{+}^{i_2}} \wedge \dots \wedge \frac{(f^{i_p})^{*} \omega}{d_{+}^{i_p}} \wedge \frac{(f^n)^{*} \omega^{s-p-l'}}{d_{+}^{(s-p-l')n}} .
  \end{split}
\end{equation*}

pour une fonction $0 \leq \phi_{l',1} \leq 1$ qui est $C^{\infty}$ à support compact dans $\Omega$. 

On a donc éliminé un  $\frac{(f^i)^{*} \omega}{d_{+}^i}$ et en recommencant ce procédé on
obtient le lemme.

\end{proof}

\begin{lemme}{\label{majoration}}

Soit $0 \leq \psi \leq 1$ une fonction $C^{\infty}$ à support compact
dans $\Omega$. Alors, on a: 

$$I=\int \psi [\xi]  \wedge \frac{(f^{i_1})^{*} \omega}{d_{+}^{i_1}} \wedge \dots \wedge \frac{(f^{i_s})^{*} \omega}{d_{+}^{i_s}}
 \leq \int [\xi] \wedge T_{+}^s + \frac{K(\xi, \psi)}{d_{+}^{i_1 }}$$

pour tous entiers $0 \leq i_1 \leq  \dots \leq i_s $.

\end{lemme}

\begin{proof}

La démonstration est du même style que la précédente: il s'agit de remplacer les $\frac{(f^{i})^{*} \omega}{d_{+}^{i}}$ par des $T_{+}$. 

On a en effet:

\begin{equation*}
\begin{split}
&I=\int \psi [\xi] \wedge \frac{(f^{i_1})^{*} \omega}{d_{+}^{i_1}} \wedge \dots \wedge \frac{(f^{i_s})^{*} \omega}{d_{+}^{i_s}}=\\
& \int \psi [\xi] \wedge T_{+} \wedge \frac{(f^{i_2})^{*} \omega}{d_{+}^{i_2}} \wedge \dots \wedge \frac{(f^{i_s})^{*} \omega}{d_{+}^{i_s}} + \int (v_{i_1}- v_{\infty}) dd^c \psi \wedge [\xi]  \wedge  \frac{(f^{i_2})^{*} \omega}{d_{+}^{i_2}} \wedge \dots \wedge \frac{(f^{i_s})^{*} \omega}{d_{+}^{i_s}} 
  \end{split}
\end{equation*}

Comme $f^r(\Omega)$ est dans $\Pp^k(\Cc) \setminus V_{+}$ pour
tout $r \geq 0$ et
que $u$ est bornée sur cet ensemble, on en déduit que $|v_{i_1} - v_{\infty}| \leq \frac{K}{d_{+}^{i_1}} $ dans $\Omega$. Cela implique que 

\begin{equation*}
\begin{split}
&I \leq \int \psi [\xi] \wedge  T_{+} \wedge \frac{(f^{i_2})^{*} \omega}{d_{+}^{i_2}} \wedge \dots \wedge \frac{(f^{i_s})^{*} \omega}{d_{+}^{i_s}}+\\
& \frac{K( \psi)}{d_{+}^{i_1}} \int \psi_1 [\xi] \wedge \omega \wedge  \frac{(f^{i_2})^{*} \omega}{d_{+}^{i_2}} \wedge \dots \wedge \frac{(f^{i_s})^{*} \omega}{d_{+}^{i_s}} 
\end{split}
\end{equation*}
 où $\psi_1$ est une fonction $C^{\infty}$ à support compact dans
$\Omega$ et comprise entre $0$ et $1$. Finalement, en utilisant le lemme précédent (quitte à renommer la constante $K(\xi, \psi)$),

$$I  \leq \int \psi [\xi] \wedge T_{+} \wedge \frac{(f^{i_2})^{*} \omega}{d_{+}^{i_2}} \wedge \dots \wedge \frac{(f^{i_s})^{*} \omega}{d_{+}^{i_s}}+ \frac{K(\xi, \psi)}{d_{+}^{i_1}} . $$

On a donc remplacé un $\frac{(f^{i})^{*} \omega}{d_{+}^{i}}$ par un $T_{+}$ et en recommençant $s$ fois on obtient le lemme.

\end{proof}

\begin{lemme}

Soit $0 \leq \psi \leq 1$ une fonction $C^{\infty}$ à support compact
dans $\Omega$.

Si on note $S_n= \frac{(f^n)_{*} (\psi
  [\xi])}{d_{+}^{sn}}$, alors

$$\| d(S_n \wedge T_{m_n,+}^p )\| \rightarrow 0 \mbox{ et } \| d d^c
  (S_n \wedge T_{m_n,+}^p )\| \rightarrow 0$$

pour toute suite $m_n$ qui tend vers l'infini et tout entier $p$
compris entre $0$ et $s-1$.

\end{lemme}

\begin{proof}

Dans le cas où $p=0$ ce lemme est démontré par N. Sibony dans le paragraphe 2.6 de
\cite{S}. Lorsque $p > 0$, il existe des versions proches de ce lemme
dans \cite{GS}. Pour le confort du lecteur, nous allons quand même en
donner une démonstration. Elle reposera comme dans \cite{S} et
\cite{GS} sur l'inégalité de Cauchy-Schwarz. Le point crucial est que
l'on pousse en avant $\xi$ qui est de dimension $s$ et que le $s$-ème
degré dynamique $d_s$ est strictement plus grand que les autres (voir
\cite{DT}).

Soit $p$ un entier compris entre $0$ et $s-1$. On va commencer par
contrôler $\| \partial(S_n \wedge T_{m_n,+}^p )\|$.

Par définition  

$$\| \partial(S_n \wedge T_{m_n,+}^p )\|= \sup_{\phi \in
  \mathcal{F}(s-p-1,s-p)} \left| \langle \partial(S_n \wedge
  T_{m_n,+}^p ) , \phi \rangle \right|, $$

où $\mathcal{F}(s-p-1,s-p)$ est l'ensemble des formes lisses $\phi$ de
bidegré $(s-p-1,s-p)$ et de norme inférieure ou égale à $1$.

{\bf Remarque:} Rappelons qu'ici $T_{m_n,+}^p$ est une forme lisse sur $\Pp^k(\Cc)
\setminus V_{+}$ et que le support de $S_n$ est aussi dans $\Pp^k(\Cc)
\setminus V_{+}$. En particulier, $S_n \wedge
  T_{m_n,+}^p$ est défini globalement par $\langle S_n \wedge
  T_{m_n,+}^p, \phi' \rangle= \langle S_n,
  T_{m_n,+}^p \wedge \phi' \rangle$, où $\phi'$ est une forme lisse de
  bidegré $(s-p,s-p)$. 

Si $\phi \in \mathcal{F}(s-p-1,s-p)$, on peut écrire

$$\phi= \sum_{i=1}^{K} \theta_i \wedge \Omega_i$$

où $K$ est une constante qui dépend seulement de $\Pp^k(\Cc)$, les
$\theta_i$ sont des $(0,1)$ formes lisses avec $\| \theta_i \| \leq 1$
et les $\Omega_i$ sont des $(s-p-1,s-p-1)$ formes lisses (strictement)
positives et de norme inférieure à $K$. Pour majorer $\left| \langle \partial(S_n \wedge
  T_{m_n,+}^p ) , \phi \rangle \right|$, il suffit donc de majorer $\left| \langle \partial(S_n \wedge
  T_{m_n,+}^p ) , \theta \wedge \Omega \rangle \right|$ où $\theta$
  est une $(0,1)$ forme lisse et $\Omega$
  est une  $(s-p-1,s-p-1)$ formes lisses (strictement)
positives avec $\| \theta\| \leq 1$ et $\| \Omega \| \leq 1$.

On a:

$$I=\left| \langle \partial(S_n \wedge  T_{m_n,+}^p ) , \theta \wedge
\Omega \rangle \right|= \left| \left\langle \frac{f_{*}^n(\partial \psi \wedge
  [\xi])}{d_{+}^{sn}} , T_{m_n,+}^p \wedge \theta \wedge \Omega \right\rangle \right|$$
 
Si $\xi' \Subset \xi$ est suffisamment proche de $\xi$, on a:

$$I=\left| \int_{f^n(\xi')} \frac{f_{*}^n(\partial \psi)}{d_{+}^{sn}}
\wedge T_{m_n,+}^p \wedge \theta \wedge \Omega \right|.$$

Pour $\alpha$ et $\beta$ des $(0,1)$ formes lisses, on définit $(\alpha, \beta):= \int_{f^n(\xi')} i \alpha
\wedge \overline{\beta} \wedge T_{m_n,+}^p  \wedge \Omega$ (voir
\cite{S}, \cite{GS} et \cite{DT}). Le fait que $(\alpha,\alpha) \geq 0$ implique
(via la démonstration de l'inégalité de Cauchy-Schwarz) que
$|(\alpha,\beta)| \leq | (\alpha,\alpha) |^{1/2}  |(\beta,\beta)
|^{1/2}$. On en déduit que

$$I \leq  \frac{1}{d_{+}^{sn}} \left| \int_{f^n(\xi')} f_{*}^n(i \partial
\psi \wedge \overline{\partial \psi})
\wedge T_{m_n,+}^p \wedge \Omega \right|^{1/2} \left| \int_{f^n(\xi')}
i \theta \wedge \overline{\theta} \wedge T_{m_n,+}^p \wedge \Omega
\right|^{1/2}.$$
Mais $i \partial \psi \wedge \overline{\partial \psi} \leq K(\psi)
\omega$ et $ \Omega \leq K \omega^{s-p-1}$ avec $K$ une constante
qui ne dépend que de $\Pp^k(\Cc)$ (car $\| \Omega \| \leq 1$). On en
déduit que

$$ \left| \int_{f^n(\xi')} f_{*}^n(i \partial
\psi \wedge \overline{\partial \psi})
\wedge T_{m_n,+}^p \wedge \Omega \right|^{1/2} \leq K'(\psi)  \left|
\int_{f^n(\xi')} f_{*}^n(\omega) \wedge T_{m_n,+}^p \wedge
\omega^{s-p-1} \right|^{1/2}$$

qui est plus petit que $K(\psi, \xi) d_{+}^{n(s-1)/2}$ d'après le premier lemme. Il reste à majorer $\left| \int_{f^n(\xi')}
i \theta \wedge \overline{\theta} \wedge T_{m_n,+}^p \wedge \Omega
\right|^{1/2}$. Mais $i \theta \wedge \overline{\theta} \wedge \Omega$
est inférieur à $K \omega^{s-p}$, et en utilisant toujours le premier lemme, on obtient:

$$\left| \int_{f^n(\xi')}
i \theta \wedge \overline{\theta} \wedge T_{m_n,+}^p \wedge \Omega
\right|^{1/2} \leq K(\psi, \xi) d_{+}^{sn/2}.$$

Finalement:

$$I \leq K(\psi, \xi) d_{+}^{-n/2}$$

(toujours quitte à renommer la constante $K(\psi, \xi)$) et cela
démontre que $\| \partial(S_n \wedge T_{m_n,+}^p )\|$ converge vers
$0$.

De là, on en déduit que $\| d(S_n \wedge T_{m_n,+}^p )\| \rightarrow
0$.

Passons à $\| dd^c (S_n \wedge T_{m_n,+}^p ) \| = \sup_{\phi \in
  \mathcal{F}(s-p-1,s-p-1)} \left| \langle dd^c(S_n \wedge
  T_{m_n,+}^p ) , \phi \rangle \right|$. Ici
  $\mathcal{F}(s-p-1,s-p-1)$ désigne l'ensemble des $(s-p-1,s-p-1)$
  formes lisses de norme inférieure ou égale à $1$.

Si $\phi$ est une telle $(s-p-1,s-p-1)$ forme, alors on peut la
décomposer en une somme $\sum_{i=0}^{K} a_i \Omega_i$ où les $a_i$
sont des nombres complexes de module plus petit que $1$ et les
$\Omega_i$ sont des $(s-p-1,s-p-1)$ formes lisses (strictement)
positives de norme inférieure ou égale à $K$ (ici $K$ ne dépend que de $\Pp^k(\Cc)$). Il suffit donc de majorer 

$$\left|  \int_{f^n(\xi')} \frac{f_{*}^n(dd^c \psi )}{d_{+}^{sn}}
\wedge T_{m_n,+}^p \wedge \Omega \right|,$$

avec $\Omega$ une $(s-p-1,s-p-1)$ forme lisse positive et $\| \Omega
\| \leq 1$ et $\xi' \Subset \xi$ suffisamment proche de $\xi$. Mais $-K(\psi) \omega \leq dd^c
\psi \leq K(\psi) \omega$ et $\Omega \leq K' \omega^{s-p-1}$ donc 

$$\left|  \int_{f^n(\xi')} \frac{f_{*}^n(dd^c \psi )}{d_{+}^{sn}}
\wedge T_{m_n,+}^p \wedge \Omega \right| \leq K'(\psi) \left|  \int_{f^n(\xi')} \frac{f_{*}^n(\omega )}{d_{+}^{sn}}
\wedge T_{m_n,+}^p \wedge \omega^{s-p-1} \right| \leq K(\psi,\xi) d_{+}^{-n}$$

par le premier lemme. Cela démontre bien que $\| dd^c (S_n \wedge
T_{m_n,+}^p ) \| \rightarrow 0 $.

\end{proof}

Les lemmes que nous venons de démontrer vont permettre de prouver
la proposition suivante.

\begin{proposition}{\label{convergence}}

Soit $0 \leq \psi \leq 1$ une fonction $C^{\infty}$ à support compact
dans $\Omega$.

Si on note $S_n= \frac{(f^n)_{*} (\psi
  [\xi])}{d_{+}^{sn}}$, alors pour tout $p$ compris entre $0$ et $s$,
  on a

$$S_n \wedge  T_{m_n,+}^p  \rightarrow c T_{-}^{k-s} \wedge T_{+}^p$$

avec $c=\int \psi [\xi] \wedge T_{+}^s$. Ici $m_n$ est une suite quelconque
d'entiers qui tend vers $+ \infty$.

\end{proposition}

\begin{proof}

Ici les wedges $S_n \wedge  T_{m_n,+}^p$ et $T_{-}^{k-s} \wedge T_{+}^p$ sont considérés de façon globale (voir la remarque au début de la preuve du lemme précédent) et pas juste définis dans un ouvert qui contient $f^n(\xi)$, sinon il n'y a pas de sens pour la convergence.

Nous allons démontrer ce résultat par récurrence sur $p$.

Supposons que $p=0$.

$S_n$ est une suite de courants positifs, de masse 

$$\int S_n \wedge \omega^{s}= \int \psi [\xi] \wedge
\frac{(f^{n})^{*} \omega^{s}}{d_{+}^{sn}} \leq K(\psi, \xi)$$

d'après le premier lemme. La suite $S_n$ admet donc des sous-suites
qui convergent. Soit $S$ une de ces limites. La masse de $S$ est égale
à $c=\int \psi [\xi] \wedge T_{+}^s$ d'après la proposition 2.6.1 de
\cite{S}.

Montrons que $S= c T_{-}^{k-s}$. Si $c=0$, $S$ est nul et le résultat est vrai. Supposons maintenant $c>0$. Le courant $S/c$ est un $(k-s,k-s)$
courant positif, fermé d'après \cite{S} (voir aussi le lemme
précédent), de masse $1$ et qui a son support dans $\overline{K_{-}}$. Mais d'après
le théorème 5.5.4 de \cite{DS5}, il est donc égal à
$T_{-}^{k-s}$. Toute les valeurs d'adhérences de $S_n$ converge donc
vers  $c T_{-}^{k-s}$ et cela démontre le cas $p=0$.

Supposons maintenant la propriété vraie au rang $p$ (avec $p < s$) et
montrons la au rang $p+1$.

Si $R$ est un courant de bidegré $(k-s+p,k-s+p)$, $v$ est une fonction
lisse et $\phi$ est une forme lisse de bidegré $(s-p-1, s-p-1)$, on a
(via des intégrations par parties):

{\bf Fait:}
$$\langle R \wedge dd^c v, \phi \rangle= \langle dd^c (vR), \phi
\rangle+ 2 \langle d R, v d^c \phi \rangle + \langle dd^c R , v \phi \rangle.$$

Considérons $R_n = S_n \wedge  T_{m_n,+}^p$. Par hypothèse de
récurrence, cette suite de courants converge vers $c T_{-}^{k-s}
\wedge T_{+}^p$.

Soit $\phi$ une forme lisse de bidegré $(s-p-1,s-p-1)$. On doit
montrer que $\langle R_n \wedge T_{m_n,+}, \phi \rangle$ converge vers
$\langle c T_{-}^{k-s}
\wedge T_{+}^{p+1} , \phi \rangle$.

On a

$$\langle R_n \wedge T_{m_n,+}, \phi \rangle=\langle R_n \wedge
\omega, \phi \rangle + \langle R_n \wedge dd^c v_{m_n}, \phi \rangle.$$

Le premier terme $\langle R_n \wedge \omega, \phi \rangle$ converge
vers $\langle c T_{-}^{k-s}
\wedge T_{+}^p \wedge \omega, \phi \rangle$. Passons au deuxième
terme. D'après le fait précédent:

$$ \langle R_n \wedge dd^c v_{m_n}, \phi \rangle = \langle
dd^c(v_{m_n}R_n), \phi \rangle + a_n+b_n$$

avec $a_n=2 \langle d R_n, v_{m_n} d^c \phi \rangle$ et $b_n= \langle dd^c R_n,
v_{m_n} \phi \rangle$.

Les courants $R_n$ vivent dans $\Pp^k(\Cc) \setminus V_{+}$. Sur cet
ensemble, les fonctions $v_{m_n}$ sont en valeur absolue bornée par
une constante $K$ indépendante de $n$. En particulier $2 v_{m_n} d^c
\phi$ et $v_{m_n} \phi$ sont des formes lisses et de norme bornée par une
constante $K(\phi)$ et le lemme précédent implique que $a_n$ et $b_n$
tendent vers $0$. Il reste à étudier $\langle dd^c(v_{m_n}R_n), \phi
\rangle= \langle v_{m_n}R_n, dd^c \phi \rangle $. Mais sur $\Pp^k(\Cc)
\setminus V_{+}$, $v_{m_n}$ converge uniformément vers $v_{\infty}$
(qui est donc continue sur cet ensemble) et
$R_n$ est une suite de courant positifs qui converge vers $c
T_{-}^{k-s} \wedge T_{+}^p$ d'où

$$\langle v_{m_n}R_n, dd^c \phi \rangle= \langle v_{\infty} R_n ,dd^c \phi
\rangle + \langle (v_{m_n} - v_{\infty}) R_n,dd^c \phi \rangle$$

converge vers $ \langle v_{\infty} c T_{-}^{k-s} \wedge T_{+}^p, dd^c
\phi \rangle$. Cela termine la démonstration de la proposition.

\end{proof}

\section{{\bf Démonstration du théorème}}

L'objectif de ce paragraphe est de démontrer que $\mu$ est l'unique mesure d'entropie maximale de $f$ dans $\Cc^k$. Nous allons suivre pour cela la méthode utilisée par E. Bedford, M. Lyubich et J. Smillie dans \cite{BLS}. Le fait de ne plus être en dimension $2$, nous obligera à changer un certain nombre d'arguments. 

Considérons $\nu$ une probabilité invariante d'entropie $\log d_+^s$. Si on fait une décomposition ergodique de la mesure $\nu$ en
$$ \nu= \int \nu_{\alpha} d \alpha ,$$
on constate que presque toutes les mesures ergodiques  $\nu_{\alpha}$ sont d'entropie $\log d_+^s$. Pour démontrer le théorème, nous sommes donc ramenés à montrer que toute probabilité invariante, ergodique et d'entropie maximale est égale à $\mu$. Dans toute la suite, nous supposerons donc $\nu$ ergodique.

Montrons tout d'abord que le support de $\nu$ est un compact de $\Cc^k$. La mesure $\nu$ peut être vue comme une probabilité de $\Pp^k$ en la prolongeant trivialement sur l'hyperplan $L_{\infty}$. Cette mesure est invariante et ne peut pas charger $V_{+}$. En effet, si $\nu(V_{+})= \alpha > 0$, alors  $\nu(\cap_{n \geq 0} f^{-n} (V_{+}))= \alpha$ car $f^{-1} (V_{+}) \subset V_{+}$. On aurait donc que $\nu(I_+)>0$ ce qui est absurde. De même, $\nu$ ne charge pas $V_{-}$. Maintenant, comme il y a un petit voisinage de $L_{\infty}$ privé de $V_{+}$ qui s'envoie par $f$ dans $V_{-}$, ce voisinage ne peut pas être chargé par $\nu$ et on en déduit que le support de $\nu$ est un compact invariant $K(\nu)$ de $\Cc^k$. En particulier, la fonction $\log d(. , L_{\infty} )$ est intégrable pour $\nu$ (ici $d$ est la distance de Fubini-Study de $\Pp^k(\Cc)$). Le Corollaire 3 de \cite{DT1} implique alors que la mesure $\nu$ est hyperbolique avec
$$\chi_1 \geq \dots \geq \chi_s \geq \frac{1}{2} \log d_+  > 0$$
et
$$ 0> - \frac{1}{2} \log d_- \geq \chi_{s+1} \geq \dots \geq \chi_k$$

où les $\chi_i$ désignent les exposants de Lyapounov de $\nu$. On a alors l'existence de variétés stables et instables locales pour la mesure $\nu$ grâce à la théorie de Pesin (toute la dynamique se passe dans un compact de $\Cc^k$).

D'après un théorème de Y. Pesin (voir \cite{Pe} et \cite{LS}), il existe une partition mesurable $f^{-1}$-invariante $\xi^{in}$ dont les fibres sont des ouverts dans les variétés instables locales et telle que
$$h_{\nu}(f)= h_{\nu}(f, \xi^{in}).$$

On va maintenant appliquer la méthode d'E. Bedford, M. Lyubich et J. Smillie (voir la proposition 3.2 de \cite{BLS}). Le point crucial pour appliquer cette méthode est la proposition suivante:

\begin{proposition}{\label{instable}}
On a
$$\int [\xi^{in}(x)] \wedge T_+^s > 0$$
pour $\nu$ presque tout point $x$.

\end{proposition}

Ici $\xi^{in}(x)$ désigne l'atome qui contient $x$ de la partition $\xi^{in}$. 

Dans \cite{BLS}, la démonstration de cette proposition passe par l'utilisation du principe du minimum pour les fonctions harmoniques (car $s=1$). Ici nous remplacerons cet argument par une démonstration qui utilise la dynamique de l'application $f$. Admettons pour l'instant cette proposition et finissons la preuve du théorème.

\subsubsection{{\bf Fin de la preuve du théorème}}

Commençons par préciser quelques notations. Pour $\nu$ presque tout point $x$,  on note $\xi^{in}(x)$ l'atome qui contient $x$ de la partition $\xi^{in}$ précédente. Par construction, $\xi^{in}(x)$ est un ouvert d'une variété instable locale (donc de dimension $s$) et il admet un prolongement $\widetilde{\xi}^{in}(x)$. En particulier,  si $U_x$ désigne un voisinage de $\xi^{in}(x)$ de sorte que $\xi^{in}(x)$ soit un sous-ensemble analytique de $U_x$, on peut définir dans $U_x$ la mesure  $T_+^s \wedge [\xi^{in}(x)]$. C'est une mesure de masse finie car les potentiels de $T_+$ sont bornés sur un voisinage de $\xi^{in}(x)$. Cette mesure dans $U_x$ peut être vue comme une mesure globale de $\Cc \Pp^k$ en posant $T_+^s \wedge  [\xi^{in}(x)](B)=T_+^s \wedge [\xi^{in}(x)](B \cap U_x)$. 

Maintenant, on peut suivre mot pour mot la preuve de la proposition 3.2 de \cite{BLS}. Il suffit de remplacer $\mu^+$  par $T_+^s $ et le degré $d$ par $d_+^s$. On obtient que pour $\nu$ presque tout point $x$, la mesure 
$$\eta_x=T_+^s \wedge [\xi^{in}(x)]/ \rho(x)$$
avec $\rho(x)= \int T_+^s \wedge [\xi^{in}(x)]> 0$ est égale à la mesure conditionnelle $\nu_x=\nu(. | \xi^{in}(x))$.

Ensuite, si on suit la preuve du théorème 3.1 de \cite{BLS}, on obtient que 
$$\frac{1}{n} \sum_{i=0}^{n-1} f^{i}_*( \nu_x) \rightarrow \nu$$
pour $\nu$ presque tout point $x$.

Comme $\nu_x$ est égale à $\eta_x$, il reste à montrer que 
$$\frac{1}{n} \sum_{i=0}^{n-1} f^{i}_*( \eta_x) \rightarrow \mu$$
et le théorème sera démontré.

Nous allons maintenant détailler cette convergence. On utilisera la définition globale de la mesure $\eta_x$ sinon quand on pousse cette mesure par $f^{i}$, on obtient quelque chose de défini dans $f^{i}(U_x)$ et ces $f^{i}(U_x)$ n'ont pas de "limite". On utilisera en particulier les résultats démontrés dans le paragraphe \ref{propriete}.

Soient $\phi$ une fonction continue et $\epsilon >0$. On considère $\psi_{\alpha}$ une fonction $C^{\infty}$ à support compact dans $U_x$ qui est comprise entre $0$ et $1$, qui vaut $1$ dans $U_x$ privé d'un $2 \alpha$-voisinage de $\partial U_x$ et $0$ dans un $\alpha$-voisinage du bord de $U_x$. On a

$$\left< \frac{1}{n} \sum_{i=0}^{n-1} f^{i}_*( \eta_x) , \phi \right>=\left<  \frac{1}{n} \sum_{i=0}^{n-1} f^{i}_*(\psi_{\alpha} \eta_x) , \phi \right> + \left< \frac{1}{n} \sum_{i=0}^{n-1} f^{i}_*( (1- \psi_{\alpha}) \eta_x) , \phi \right>.$$
Le second terme de la somme est égal à $\left<  (1- \psi_{\alpha}) \eta_x , \frac{1}{n} \sum_{i=0}^{n-1} \phi \circ f^{i}  \right>$, qui en valeur absolue est inférieur à $ \| \phi \| \eta_x( 2 \alpha- \mbox{voisinage de } U_x)$ qui est très petit devant $\epsilon$ si on prend $\alpha$ petit (les bords de $\xi^{in}$ peuvent être pris génériques par rapport à $T_+^s$ si on veut).

Passons au premier terme de la somme précédente.

Il est égal à

$$\frac{1}{\rho(x)}  \frac{1}{n} \sum_{i=0}^{n-1}   \left< \frac{(f^{i})^* (T_+^s)}{d_+^{si}} \wedge [\xi^{in}(x)]  \psi_{\alpha} , \phi \circ f^{i} \right>$$

car $f^* T_+^s=d_+^s T_+^s$. De plus, on a:

{\bf Fait:}

$$\  \left< \frac{(f^{i})^* (T_{m,+}^s)}{d_+^{si}} \wedge [\xi^{in}(x)]  \psi_{\alpha} , \phi \circ f^{i} \right>$$

converge vers

$$ \left< \frac{(f^{i})^* (T_{+}^s)}{d_+^{si}} \wedge [\xi^{in}(x)]  \psi_{\alpha} , \phi \circ f^{i} \right>$$

quand $m$ tend vers l'infini.

La démonstration de ce fait est classique et passe si on veut par l'écriture $T_{m,+}= \omega + dd^c v_m$  et $T_{+}= \omega + dd^c v_{\infty}$ (voir aussi le paragraphe \ref{propriete}).

Maintenant, si on prend des $m_i$ suffisamment grands par rapport à $i$ et $\alpha$, on a que la différence entre
$$\left< \frac{1}{n} \sum_{i=0}^{n-1} f^{i}_*( \eta_x) , \phi \right>$$ 
et
$$\frac{1}{\rho(x)}  \frac{1}{n} \sum_{i=0}^{n-1}   \left< \frac{(f^{i})^* (T_{m_i,+}^s)}{d_+^{si}} \wedge [\xi^{in}(x)]  \psi_{\alpha} , \phi \circ f^{i} \right>$$
est petite devant $\epsilon$.

Enfin, ce dernier terme est égal à:

$$\frac{1}{\rho(x)}  \frac{1}{n} \sum_{i=0}^{n-1}   \left< T_{m_i,+}^s \wedge \frac{f^{i}_*(\psi_{\alpha} [\xi^{in}(x)])}{d_+^{si}} , \phi  \right>$$
 
qui converge vers

$$\left( \frac{\int  \psi_{\alpha} [\xi^{in}(x)] \wedge T_+^s }{\rho(x)} \right)  \left< \mu , \phi \right>$$

quand $n$ tend vers l'infini grâce à la proposition \ref{convergence}.

Ce dernier est aussi proche que l'on veut de $\left< \mu , \phi \right>$ si on prend $\alpha$ petit.

On a donc bien montré que 

$$\frac{1}{n} \sum_{i=0}^{n-1} f^{i}_*( \eta_x) \rightarrow \mu$$

et le théorème est démontré.

Il reste à prouver la proposition \ref{instable}: ce sera l'objet du paragraphe suivant.

\subsubsection{{\bf Démonstration de la proposition}}

Dans ce paragraphe nous allons montrer que $\int [\xi^{in}(x)] \wedge T_+^s > 0$ pour $\nu$ presque tout point $x$.

Ici les arguments ne sont pas les mêmes que ceux utilisés par E. Bedford, M. Lyubich et J. Smillie. En effet, dans leur situation $s=1$. Le fait que $\int [\xi^{in}(x)] \wedge T_+ = 0$ implique que le potentiel de $T_+$ est harmonique sur $\xi^{in}(x)$. Ils en déduisent alors une contradiction en utilisant le principe du minimum. Ce raisonnement ne peut pas être fait en dimension supérieure. Pour le remplacer, nous allons utiliser des arguments dynamiques qui  vont utiliser une idée de S. Newhouse (voir \cite{Ne}). Une autre approche possible m'a été signalée par R. Dujardin pour démontrer ce point là. Elle s'appuie sur la démonstration de la proposition 5.1 de \cite{Du}, en y injectant le théorème de Bézout pour les courants dans $\Pp^k(\Cc)$. Nous espérons que notre approche "plus locale" pourra être étendue à certaines applications holomorphes ou méromorphes dans les variétés Kählériennes compactes.

Soit

$$M:= \{ x \mbox{, } \int [\xi^{in}(x)] \wedge T_+^s = 0 \}.$$

On va montrer que $M$ est un ensemble négligeable pour $\nu$. Soit $\epsilon >0$. La mesure $\nu$ est hyperbolique avec $s$ exposants strictement positifs et $k-s$ strictement négatifs: pour $\nu$ presque tout point on a donc l'existence d'une variété stable locale $\xi^{st}(x)$ de dimension $k-s$  et d'une variété instable locale $\xi^{in}(x)$ de dimension $s$. Ces variétés sont des graphes de fonctions $\phi^{st}_x:   B^{st}(x, r^{st}(x))  \rightarrow E^{in}(x)$ (respectivement $\phi^{in}_x:B^{in}(x, r^{in}(x)) \rightarrow E^{st}(x)$) (voir le paragraphe 2.4 de \cite{BLS} pour les notations et des rappels sur la théorie de Pesin). Si on choisit $\alpha_0 >0$ suffisamment petit et $A$ suffisamment grand, on a que sur un ensemble $\Lambda_{\alpha_0}$ de $x$ de mesure supérieure à $1 - \epsilon/4$ pour $\nu$: 

- $r^{in}(x)$, $r^{st}(x)$ et l'angle entre $E^{st}(x)$ et $E^{in}(x)$ sont supérieurs à $\alpha_0 >0$.

- $\| \phi^{st}_x \|_{C^2} \leq A$ sur $B^{st}(x, r^{st}(x))$ et de même  $\| \phi^{in}_x \|_{C^2} \leq A$ sur $B^{in}(x, r^{in}(x))$.

- $x \rightarrow E^{st}(x)$ et $x \rightarrow E^{in}(x)$ sont continues sur $\Lambda_{\alpha_0}$ (théorème de Lusin).

- pour $y$ dans $\xi^{st}(x)$ et  $n \geq 0$, on a $dist(f^n(x), f^n(y)) \leq A  e^{-(\lambda - \gamma)n}$ (avec $\lambda= \min \{ | \chi_i | \mbox{, } \chi_i < 0 \}$ et $\gamma$ petit par rapport à $\lambda$).

 Maintenant, par le théorème de Brin et Katok (voir \cite{BK}), on a 
 $$h_{\nu}(f)= \log d_+^s= \lim_{\delta \rightarrow 0} \liminf_n - \frac{1}{n} \log \nu(B_n(x, \delta))$$
 
 pour $\nu$ presque tout point $x$. Donc, si on note
 
 $$\Lambda_{\delta , n}= \{ x \mbox{, } \nu(B_n(x , \delta)) \leq d_+^{-sn+ \gamma n} \} ,$$
 
 si on choisit $\delta$ assez petit (dans la suite on prendra aussi $\delta$ très petit devant $\alpha_0$ et $A$), on a
 
 $$1 - \frac{\epsilon}{8} \leq \nu \left( \left\{ x \mbox{, } \liminf_n - \frac{1}{n} \log \nu(B_n(x, \delta)) \geq  (s- \gamma/2) \log d_+ \right\} \right) \leq \nu( \cup_{n_0} \cap_{n \geq n_0}   \Lambda_{\delta , n} ).$$
 
 En particulier, si on prend $n_0$ grand, on a $\nu(  \cap_{n \geq n_0}   \Lambda_{\delta , n} ) \geq 1 - \frac{\epsilon}{4}$. Dans la suite, on notera
 $$\Lambda= \Lambda_{\alpha_0} \cap \left( \cap_{n \geq n_0}   \Lambda_{\delta , n} \right)$$
 qui est de mesure supérieure à $1 - \frac{\epsilon}{2}$ pour $\nu$.
 
 Il suffit maintenant de montrer que $M \cap \Lambda$ est inclus dans un ensemble de mesure inférieure à $\epsilon/2$.
 
Soit $x_1 , x_2 , \dots , x_N$ un ensemble $(n, \delta)$ séparé de cardinal maximal dans $M \cap \Lambda$. On a

$$ M \cap \Lambda \subset \cup_{i=1}^{N} B_n(x_i, \delta)$$

et pour $n \geq n_0$

$$\nu(\cup_{i=1}^{N} B_n(x_i, \delta)) \leq N d_+^{-sn+ \gamma n}. $$

Si on montre que $N$ est inférieur à $d_+^{sn -2 \gamma n} $, alors on aura bien montré que $M \cap \Lambda$ est inclus dans un ensemble de mesure plus petite que $\epsilon/2$. Il nous reste donc à majorer $N$.

On considère un découpage d'un voisinage du support de $\nu$ en cubes disjoints de taille $r=r(\delta, \alpha_0, A)$ que l'on va préciser au fur et à mesure. Fixons un de ces cubes $C(r)$. Si $r$ est suffisamment petit (par rapport à $\alpha_0$) on a que pour $x$ et $y$ dans $\Lambda$ qui sont dans le même cube $C(r)$, alors les directions $E^{st}(x)$ et $E^{st}(y)$ sont très proches (disons à une distance plus petite que $\alpha_0 /1000$). De même pour $E^{in}(x)$ et $E^{in}(y)$ (cela provient de la continuité de $x \rightarrow E^{st}(x)$ et $x \rightarrow E^{in}(x)$ sur $\Lambda_{\alpha_0}$ ). Par conséquent, quand on prend deux points $x$ et $y$ dans $\Lambda$ qui sont dans le même cube $C(r)$, l'espace affine  $E^{in}(x)$ rencontre $E^{st}(y)$ en un point qui se trouve dans un cube de taille $C(\alpha_0)r$ (centré comme $C(r)$) et de même $E^{in}(y)$ rencontre $E^{st}(x)$. Ici $C(\alpha_0)$ dépend de l'angle entre les espaces stables et instables et celui-ci est minoré par essentiellement $\alpha_0$.

Ces points d'intersections vont impliquer que d'une part $\xi^{st}(x)$ rencontre $\xi^{in}(y)$ dans le cube de taille $2C(\alpha_0)r$ et d'autre part $\xi^{in}(x)$ rencontre $\xi^{st}(y)$ dans ce même cube. En effet soit par exemple $z$ le point d'intersection de $E^{in}(x)$ avec $E^{st}(y)$. On considère le polydisque $P=B_1 \times B_2$ de taille $r$ centré en $z$ construit à partir de ces deux espaces ($B_1$ est la boule de centre $z$ et de rayon $r$ dans $E^{in}(x)$ et $B_2$ celle de centre $z$ et de rayon $r$ dans  $E^{st}(y)$). Maintenant,  dans le cube de taille $2C(\alpha_0)r$ (qui contient $P$ si on veut) la différentielle de $\phi^{in}_x $ est majorée par $2 \sqrt{k} A C(\alpha_0)r$ car la différentielle seconde est majorée par $A$ et que la variété instable est tangente à $E^{in}(x)$  en $x$. De même pour la variété stable en $y$. En particulier, dans le cube $2C(\alpha_0)r$, la distance entre le graphe de $\phi^{in}_x $ et $E^{in}(x)$  est plus petite que $4 k A (C(\alpha_0)r)^2$ (de même pour l'autre). Comme  $A C(\alpha_0)r$ peut être pris petit (si on prend $r=r(\delta, \alpha_0, A)$), on en déduit que la variété instable en $x$, $\xi^{in}(x)$ rencontre la variété stable $\xi^{st}(y)$ de $y$ en (au moins) un point (voir par exemple la proposition S.3.7 de \cite{KH}). De la même façon,  $\xi^{st}(x)$ rencontre $\xi^{in}(y)$. Les points d'intersections se trouvent dans le cube $C$ de taille $2C(\alpha_0)r$. Maintenant, on peut choisir $r$ petit par rapport à $\delta$ de sorte que le diamètre de $f^n(\xi^{st}(x) \cap C)$ soit très petit devant $\delta/4$ pour tout $n \geq 0$ et $x$ dans le cube $C(r) \cap \Lambda$. En effet, si $n$ est grand cela provient directement du contrôle $dist(f^n(x), f^n(y)) \leq A  e^{-(\lambda - \gamma)n}$ pour $y$ dans $\xi^{st}(x)$. Pour les autres $n$ il suffit de réduire la taille de la boîte (i.e. prendre $r$ petit) et d'utiliser que $f$ est holomorphe et le théorème des accroissement finis.

Si on reprend les $N$ points $x_1 , x_2 , \dots , x_N$ de $M \cap \Lambda$ que l'on avait, on peut trouver un cube de taille $r$, $C(r)$ qui contient au moins $Nr^{2k}$  points $x_i$. Fixons un point $x$ dans $C(r) \cap M \cap \Lambda$ indépendant de $n$ (i.e. on fixe un $x$ dans tous les cubes qui contiennent un point de $M \cap \Lambda$). Sa variété instable $\xi^{in}(x)$ intersectée avec le cube $C$ de taille $2C(\alpha_0)r$ qui est centré comme $C(r)$,  rencontre les variétés stables $\xi^{st}(x_j) \cap C$ des $x_j$ qui sont dans $C(r)$. Notons $z_1 , z_2 , \dots , z_{N_0}$ ces points d'intersection avec la convention de ne prendre qu'un seul point d'intersection entre $\xi^{in}(x) \cap C$ et $\xi^{st}(x_j) \cap C$ si jamais il y en a plusieurs. On a $N_0 \geq N r^{2k}$. 

C'est ce nombre $N_0$ que nous allons maintenant majorer en utilisant un argument d'entropie. En effet, remarquons tout d'abord (comme dans \cite{Ne}) que les points $z_1 , z_2 , \dots , z_{N_0}$ sont $(n, \delta/2)$ séparés. En effet prenons par exemple $z_1$ et $z_2$. On a $z_1 \in \xi^{st}(x_j) \cap C$ et  $z_2 \in \xi^{st}(x_l) \cap C$ avec $j \neq l$. Comme $x_j$ et $x_l$ sont $(n, \delta)$ séparés, il existe un entier $p$ entre $0$ et $n-1$ tel que $d(f^p(x_j), f^p(x_l))   \geq \delta$. Mais comme les diamètres de $f^p(\xi^{st}(x_j) \cap C)$ et de $f^p(\xi^{st}(x_l) \cap C)$ sont inférieurs à $\delta/4$ on en déduit bien que la distance entre $f^p(z_1)$ et $f^p(z_2)$ est supérieure à $\delta/2$. En conclusion, si on sait majorer le cardinal d'un ensemble $(n , \delta/2)$ séparé dans la variété instable  $\xi^{in}(x) \cap C$ par $d_+^{sn -2 \gamma n} $, alors $N$ sera majoré par $r^{-2k} d_+^{sn -2 \gamma n} $ et la proposition sera démontrée.

Il nous reste donc à majorer le cardinal d'un ensemble $(n , \delta/2)$ séparé dans la variété instable  $\xi^{in}(x) \cap C$ (où $x$ est dans $M \cap \Lambda$). Cela va se faire en deux étapes. Dans la première nous allons majorer l'entropie par un volume comme dans \cite{Gr}. Ensuite ce volume sera majoré en utilisant que $\int [\xi^{in}(x)] \wedge T_+^s = 0 $ et la convergence vers $T_+^s $. C'est une idée utilisée dans \cite{DT0}, pour montrer que pour un endomorphisme holomorphe de $\Pp^2(\Cc)$ de degré $d$, l'entropie topologique en dehors du support de la mesure de Green est majorée par $\log(d)$.

$$ $$

{\bf Majoration de l'entropie par un volume}

On va utiliser ici l'argument de Gromov (voir \cite{Gr}). Notons $C_{\delta}$ le cube de taille $2 C(\alpha_0) r+ \delta$ qui est centré comme $C$ et $\xi_{\delta}^{in}(x)$ l'intersection de $\xi^{in}(x)$ avec $C_{\delta}$. Comme $\delta$ est petit devant $\alpha_0$ et $A$, la distance entre le bord de $\xi^{in}(x)$ et le cube $C_{\delta}$ est plus grande que $\sqrt{k} \delta$.

On considère le multigraphe $\Gamma_n = \{ (y, f(y), \dots , f^{n-1}(y) ) \mbox{, } y \in \xi_{\delta}^{in}(x)  \}$. L'ensemble des points $(n, \delta/2)$ séparés $z_1 , z_2 , \dots , z_{N_0}$ induisent, via leurs $n$-orbites, un ensemble $F$ de $\Gamma_n$ qui est $\delta/2$ séparé dans $(\Cc \Pp^k)^n$ pour la métrique produit. Cela signifie que l'on a $N_0$ boules disjointes $B(a, \delta/4)$ avec $a \in F$ dans $(\Cc \Pp^k)^n$. Par le théorème de Lelong, le volume de $\Gamma_n \cap B(a, \delta/4)$ est minoré par une constante $c(\delta)$ car le bord de $\Gamma_n$ est en dehors de la boule $B(a, \delta/4)$. En effet, comme les points $z_1 , z_2 , \dots , z_{N_0}$ sont dans $C$, la distance entre ces points et le bord de $\xi_{\delta}^{in}(x) $ est supérieure à $\delta$. On vient donc de montrer que le volume de $\Gamma_n$ est supérieur à $c(\delta) N_0$. Il reste donc à majorer ce volume par $d_+^{sn -2 \gamma n} $.

$$ $$

{\bf Majoration du volume de $\Gamma_n$}

Dans $(\Cc \Pp^k)^n$, on notera $\Pi_i$ les projections sur les facteurs du produit et on munit $(\Cc \Pp^k)^n$ de la forme de Kähler $\omega_n= \sum_{i=1}^n \Pi^*_i \omega$. Le volume de $\Gamma_n$ est alors égal à
$$\int_{\Gamma_n} (\omega_n)^s = \sum_{0 \leq n_1 , \dots , n_s \leq n-1} \int_{\xi_{\delta}^{in}(x) } (f^{n_1})^* \omega \wedge \dots \wedge (f^{n_s})^* \omega.$$
Il s'agit donc de majorer les termes $\int_{\xi_{\delta}^{in}(x) } (f^{n_1})^* \omega \wedge \dots \wedge (f^{n_s})^* \omega$. Soit $\psi$ une fonction $C^{\infty}$ à support compact comprise entre $0$ et $1$, qui vaut $1$ dans le cube $C_{\delta}$ et $0$ en dehors du cube de taille $2 C(\alpha_0) r+  2 \delta$ centré comme $C$ et  $C_{\delta}$. Il s'agit donc de majorer les termes
$$\int \psi [\xi^{in}(x)]  \wedge (f^{n_1})^* \omega \wedge \dots \wedge (f^{n_s})^* \omega.$$

On va considérer pour cela deux cas: soit tous les $n_i$ sont supérieurs ou égaux à $n/2$, soit il existe un des $n_i$ inférieur à $n/2$.

$$ $$

{\bf $1^{er}$ Cas:} tous les $n_i$ sont supérieurs à $n/2$:

En utilisant le lemme \ref{majoration}, on obtient:

$$\int \psi [\xi^{in}(x)]  \wedge (f^{n_1})^* \omega \wedge \dots \wedge (f^{n_s})^* \omega \leq d_{+}^{n_1 + \dots + n_s} \left( 0 + \frac{K(\xi^{in}(x), \psi)}{d^{n/2}}  \right)$$

car $x$ est dans $M$.

On en déduit que

$$\int \psi [\xi^{in}(x)]  \wedge (f^{n_1})^* \omega \wedge \dots \wedge (f^{n_s})^* \omega \leq K(\xi^{in}(x), \psi) d_{+}^{ns - n/2}. $$

{\bf $2^{eme}$ Cas:} un des $n_i$ est inférieur à $n/2$:

Toujours par le lemme \ref{majoration}, on a:

$$\int \psi [\xi^{in}(x)]  \wedge (f^{n_1})^* \omega \wedge \dots \wedge (f^{n_s})^* \omega \leq K(\xi^{in}(x), \psi) d_{+}^{n_1 + \dots + n_s} $$

qui est majoré par $K(\xi^{in}(x), \psi) d_{+}^{ns - n/2}$ car un des $n_i$ est plus petit que $n/2$.

$$ $$

Cela conclut la démonstration de la proposition et donc celle du théorème.

\bigskip

Henry de Thélin

Université Paris-Sud (Paris 11)

Mathématique, Bât. 425

91405 Orsay

France

\end{document}